\documentclass{amsproc}
\usepackage{amssymb}

\usepackage{graphicx}
\usepackage{amscd}

\theoremstyle{plain}

\newtheorem {Theorem}   {Theorem} 

\numberwithin{Theorem}{section}

\newtheorem {Lemma}[Theorem]    {Lemma}         
  
\theoremstyle{Definition}
\newtheorem{Definition}[Theorem]{Definition}
\theoremstyle{remark}

\begin{document}
\title[Hagopian's theorem and exponents]{A generalization of Hagopian's theorem and exponents}
\author{Alex Clark}
\address{University of North Texas\\Denton, TX 76203-1430}
\email{alexc@unt.edu}
\thanks{This work was funded in part by a faculty
research grant from the University of North Texas.}
\date{March 27, 2000}
\subjclass{Primary 54F15, 37B45; Secondary 54H11}
\keywords{homogeneous, continuum, solenoid, topological group}

\begin{abstract}
We generalize Hagopian's theorem characterizing solenoids to higher dimensions 
by showing that any homogeneous continuum admitting a fiber bundle projection onto a torus with
totally disconnected fibers admits a compatible topological group structure. And then the higher
dimensional exponent group is introduced.
\end{abstract}

\maketitle

\section{Introduction}

A space $X$ is a \emph{continuum} if it is compact, connected and
metrizable, and $X$ is \emph{homogeneous} if given any $x,y\in X$ there is a
homeomorphism $h:X\rightarrow X$ with $h\left( x\right) =y$. We show that if 
$X$ is a homogeneous continuum which (for some countable $\kappa \geq 1$)
admits a fiber bundle projection $p:X\rightarrow \mathbf{T}^{\kappa }=\left( 
\mathbf{R}/\mathbf{Z}\right) ^{\kappa }$ with totally disconnected fibers,
then $X$ also admits a compatible abelian topological group structure. This
generalizes a weakened version of the following theorem of Hagopian: If
every subcontinuum of the homogeneous continuum $X$ is an arc, then $X$ is
homeomorphic to a solenoid (which includes the circle as a possibility) \cite
{H}. For $\kappa =1$, the theorem currently under consideration follows from
Hagopian's theorem and applies to all compact one--dimensional minimal sets
of flows (see \cite{AM}). In higher dimensions, spaces admitting such a
fiber bundle structure arise naturally in two settings: as the minimal sets
of foliations and as the limit sets of discrete dynamical systems (see \cite
{W}).

A good example of an $S^{1}$ fiber bundle which is not homogeneous is the
minimal set of a Denjoy flow on a torus: the path components are not ``evenly
spaced,'' and so there is no Effros homeomorphism for sufficiently small
numbers. If one follows path components far enough along, eventually they
will spread apart. Our result shows that this situation persists in higher
dimensions.

We then generalize the exponent group introduced in \cite{Cl} to higher
dimensions. In what follows we use the terminology and results of \cite{S},
and $\pi ^{\kappa }:\mathbf{R}^{\kappa }\rightarrow \mathbf{T}^{\kappa }$
denotes the standard fibration with unique path lifting $\left\langle
t_{i}\right\rangle _{i=1}^{\kappa }\longmapsto \left\langle t_{i}\,\text{mod 
}1\right\rangle _{i=1}^{\kappa }$ and $d$ is a metric for $X.$

\section{Homogeneity}

\begin{Theorem}
If $X$ is a homogeneous continuum and if $X$ admits a fiber bundle
projection $p:X\rightarrow \mathbf{T}^{\kappa }$ with totally disconnected
fibers, then $X$ admits an abelian topological group structure.
\end{Theorem}

\noindent \textbf{Proof}: Suppose $p:X\rightarrow \mathbf{T}^{\kappa }$ is a fiber
bundle projection with totally disconnected fiber $F$ for some $\kappa >1$.
Then $p$ is a fibration with unique path lifting (see \cite{S} 2.2.5,
2.7.14) and so the natural group action 
\[
\alpha :\mathbf{R}^{\kappa }\times \mathbf{T}^{\kappa }\rightarrow \mathbf{T}
^{\kappa }
\]
given by 
\[
\alpha \left( s,x\right) =\pi ^{\kappa }\left( s\right) +x
\]
lifts uniquely to an action $\tilde{\alpha}$ on $X$ 
\[
\begin{array}{lll}
\mathbf{R}^{\kappa }\times X & \stackrel{\tilde{\alpha}}{\longrightarrow } & 
X \\ 
id\times p\downarrow  &  & \downarrow p \\ 
\mathbf{R}^{\kappa }\times \mathbf{T}^{\kappa } & \stackrel{\alpha }{
\longrightarrow } & \mathbf{T}^{\kappa }
\end{array}
.
\]
For, given any $\left( s,x\right) \in \mathbf{R}^{\kappa }\times X$, we have
the interval $I=\left[ t\cdot s,\left( 1-t\right) \cdot s\right] $ ($t\in 
\left[ 0,1\right] $) and the following commutative diagram 
\[
\begin{array}{lllll}
\left\{ \mathbf{0}\right\} \times X &  & \stackrel{``id"}{\longrightarrow }
& \stackrel{}{} & X \\ 
\downarrow \cap  &  &  &  & \downarrow p \\ 
I\times X & \stackrel{id\times p}{\longrightarrow } & I\times \mathbf{T}
^{\kappa } & \stackrel{\alpha }{\longrightarrow } & \mathbf{T}^{\kappa }
\end{array}
,
\]
and so the fibration property allows us to define $\tilde{\alpha}
\upharpoonleft I\times X$ as needed. To see that $\tilde{\alpha}$ is
continuous, we can replace $I$ in the diagram with increasing open connected
subsets of $\mathbf{R}^{\kappa }$ whose union is all of $\mathbf{R}^{\kappa }
$.

Then for a given point $e\in X$, 
\[
\pi ^{\kappa }=p\circ \tilde{\alpha}\upharpoonleft \mathbf{R}^{\kappa
}\times \left\{ e\right\}
\]
is a fibration and $p$ has unique path lifting, and so it follows from lemma
2.1 of \cite{C} that $\tilde{\alpha}\upharpoonleft \mathbf{R}^{\kappa
}\times \left\{ e\right\} $ is a fibration. Hence, the $\tilde{\alpha}-$
trajectory of $e$ is the path component of $e$ in $X$, \cite{S} 2.3.1. We
now proceed to show that the path components of $X$ are dense in $X$.

The proof is similar to the proof of a related fact for matchbox manifolds
given in \cite{AHO}. Each point of $X$ is contained in an open set
homeomorphic to $V\times F$ for some open subset $V$ of $\mathbf{T}^{\kappa
} $. Since $F$ is compact and totally disconnected, there is a basis $
\mathcal{B}$ for $X$ consisting of sets homeomorphic to $U\times Z$ with $U$
open and connected in $\mathbf{T}^{\kappa }$ and $Z$ a closed and open
subset of $F.$ Let $B\approx U\times Z$ be any element of $\mathcal{B}$. We
wish to show that $PC\left( B\right) $ (the union of all the path components
of $X$ which meet $B$) is both open and closed in $X$. For $t\in \mathbf{R}
^{\kappa }-\left\{ \mathbf{0}\right\} $, $\tilde{\alpha}$ induces the
non-singular flow 
\[
\tilde{\alpha}\left[ t\right] :\mathbf{R\times }X\rightarrow X\,;\text{ } 
\tilde{\alpha}\left[ t\right] \left( r,x\right) =\tilde{\alpha}\left(
tr,x\right) \text{ (and similarly for }\alpha \text{).}
\]
For $x\in PC\left( B\right) $ we may then choose $s\in \mathbf{R}^{\kappa }$
so that $\tilde{\alpha}\left( s,x\right) =y\in B$ and so that the
corresponding linear flow $\alpha \left[ s\right] $ (and hence $\tilde{
\alpha }\left[ s\right] $) is aperiodic, which is possible since $\kappa >1$
and since the collection of $t\in \mathbf{R}^{\kappa }-\left\{ \mathbf{0}
\right\} $ for which $\alpha \left[ t\right] $ is aperiodic is dense. By a
theorem of Bebutov (see \cite{NS} V.2.15) we may construct a flow box
joining $y$ and $x$ which contains both points in its interior. Hence, all
points in a neighborhood of $x$ are also in $PC\left( B\right) $,
demonstrating that it is open.

Suppose then that $x\in \overline{PC\left( B\right) }$ with $\left(
x_{n}\right) _{n}\rightarrow x$ and $\left\{ x_{n}\right\} \subset PC\left(
B\right) .$ The path component of $x_{n}$ meets $B$ in a point $b_{n}\approx
\left( u_{n},z_{n}\right) $ and since $U$ is path connected, we may assume
that $u_{n}=u$ for some $u\in U$ and for all $n$. The compactness of $Z$
allows us to find a subsequence of $\left( z_{n}\right) $ converging to some 
$z\in Z$ with $b\approx \left( u,z\right) \in B$. Then there is some $
\varepsilon >0$ so that $B_{d}\left( b,\varepsilon \right) \subset B$.
Applying the Effros theorem \cite{E},\cite{U} to $H\left( X\right) $, the
group of homeomorphisms of $X$ in the $\sup $ metric, there is a $\delta >0$
so that for any $v,w\in X$ with $d\left( v,w\right) <\delta $ there is an $
h^{\left[ v,w\right] }\in H\left( X\right) $ with $h^{\left[ v,w\right]
}\left( v\right) =w$ and $d\left( p,h^{\left[ v,w\right] }\left( p\right)
\right) <\varepsilon /2$ for all $p\in X$ (denoted: $\delta \stackrel{\text{
Eff}}{\sim }\varepsilon /2$). We now choose $n$ so that $d\left(
b_{n},b\right) <\varepsilon /2$ and $d\left( x_{n},x\right) <\delta $ and a
corresponding homeomorphism $h^{\left[ x_{n},x\right] }$. Then $h^{\left[
x_{n},x\right] }\left( b_{n}\right) \in B_{d}\left( b,\varepsilon \right)
\subset B$ is a point in the path component of $x$, demonstrating that $x\in
PC\left( B\right) ,$ and so $PC\left( B\right) =X$. Since $B$ was any such
basis element, each path component of $X$ is dense.

With $\tilde{\alpha}_{x}:\mathbf{R}^{\kappa }\rightarrow X$ denoting the $
\tilde{\alpha}-$orbit of $x\in X$ and giving 
\[
A\stackrel{\text{def}}{=}\left\{ \tilde{\alpha}_{x}\mid x\in X\right\}
\]
the topology it inherits from the collection of all maps $\mathbf{R}^{\kappa
}\rightarrow X$ in the $\sup $ metric, we proceed to show that $
h:X\rightarrow A;$ $h\left( x\right) \stackrel{\text{def}}{=}\tilde{\alpha}
_{x}$ is a homeomorphism. Since $h$ is clearly a bijection and $X$ is
compact, it suffices to show that $h$ is continuous. Given $\left(
x_{n}\right) _{n}\rightarrow x$ in $X$ and $\varepsilon >0$ we need to find
an $N$ so that $\sup\limits_{s\in \mathbf{R}^{\kappa }}d\left( \tilde{\alpha}
_{x_{n}}\left( s\right) ,\tilde{\alpha}_{x}\left( s\right) \right)
<\varepsilon $ for all $n\geq N$.

Let $d_{\kappa }$ be a translation invariant metric for $\mathbf{T}^{\kappa }
$. First we find a Lebesgue number $\lambda >0$ for a covering $\mathcal{O}$
of $\mathbf{T}^{\kappa }$ by open sets $V$ satisfying 
\[
p^{-1}\left( V\right) \approx V\times F\text{ }
\]
and $p\left( y\right) =v$ for $y\approx \left( v,f\right) .$ \ Since $X$ and 
$\mathbf{T}^{\kappa }$ are compact, we may find a connected neighborhood $U$
of $\mathbf{0}\in \mathbf{R}^{\kappa }$ so that 
\[
\sup_{y\in \mathbf{T}^{\kappa }}\left\{ \text{diam}\alpha \left( U\times
\left\{ y\right\} \right) \right\} <\lambda \text{ and }\sup_{y\in X}\left\{ 
\text{diam}\tilde{\alpha}\left( U\times \left\{ y\right\} \right) \right\} <
\frac{\varepsilon }{3}
\]
and so that $\alpha \left( U\times \left\{ e\right\} \right) $ is the $\eta -
$neighborhood of $e$ in $\mathbf{T}^{\kappa }$ for some $\eta >0$. The
translation invariance of $d_{\kappa }$ yields that $\alpha \left( U\times
\left\{ y\right\} \right) $ is the $\eta -$neighborhood of $y\in \mathbf{T}
^{\kappa }$. Since $p$ is uniformly continuous, there is a $\tau >0$ so that 
\[
d\left( y,y^{\prime }\right) <\tau \Rightarrow d_{\kappa }\left( p\left(
y\right) ,p\left( y^{\prime }\right) \right) <\eta .
\]
Next we find 
\[
0<\delta \stackrel{\text{Eff}}{\sim }\min \left\{ \frac{1}{2}\varepsilon
,\tau \right\} \text{ and }N\text{ so that }\left\{ x_{n}\right\} _{n\geq
N}\subset B_{d}\left( x,\delta \right) .
\]
Given $n\geq N$ there is an $h^{\left[ x_{n},x\right] }\in H\left( X\right) $
within $\min \left\{ \frac{1}{2}\varepsilon ,\tau \right\} $ of $id_{X}$.
And so 
\[
p\left( x\right) =p\left( h^{\left[ x_{n},x\right] }\left( x_{n}\right)
\right) \in \alpha \left( U\times \left\{ p\left( x_{n}\right) \right\}
\right) 
\]
since $d\left( x_{n},x\right) <\tau .$ For $s\in \mathbf{R}^{\kappa }$, the
translation invariance of $d_{\kappa }$ yields that $d_{\kappa }\left(
\alpha \left( s,p\left( x\right) \right) ,\alpha \left( s,p\left(
x_{n}\right) \right) \right) <\eta $, while the choice of $h^{\left[ x_{n},x
\right] }$ and the equality $\alpha \left( s,p\left( x_{n}\right) \right)
=p\circ \tilde{\alpha}\left( s,x_{n}\right) $ yields that 
\[
d_{\kappa }\left( \alpha \left( s,p\left( x_{n}\right) \right) ,p\circ h^{
\left[ x_{n},x\right] }\circ \tilde{\alpha}\left( s,\left( x_{n}\right)
\right) \right) <\eta .
\]
Combining this with the equality $\alpha \left( s,p\left( x\right) \right)
=p\circ \tilde{\alpha}\left( s,x\right) ,$ we obtain 
\[
\left( \ast \right) \left\{ p\circ \tilde{\alpha}\left( s,x\right) ,p\circ
h^{\left[ x_{n},x\right] }\circ \tilde{\alpha}\left( s,\left( x_{n}\right)
\right) \right\} \subset \alpha \left( U\times \left\{ \alpha \left(
s,p\left( x_{n}\right) \right) \right\} \right) \text{ for all }s\in \mathbf{
R}^{\kappa }.
\]
For $s\in \mathbf{R}^{\kappa }$ let 
\[
U\left( s\right) \stackrel{\text{def}}{=}p^{-1}\left( \alpha \left( U\times
\left\{ \alpha \left( s,p\left( x_{n}\right) \right) \right\} \right)
\right) .
\]
By its size, we know that $\alpha \left( U\times \left\{ \alpha \left(
s,p\left( x_{n}\right) \right) \right\} \right) $ fits inside some $V\in 
\mathcal{O}$, and so 
\[
U\left( s\right) \approx \alpha \left( U\times \left\{ \alpha \left(
s,p\left( x_{n}\right) \right) \right\} \right) \times F
\]
and for $y\approx \left( v,f\right) $, $p\left( y\right) =v$. Then by $
\left( \ast \right) $ we have 
\[
\left\{ \tilde{\alpha}\left( s,x\right) ,h^{\left[ x_{n},x\right] }\circ 
\tilde{\alpha}\left( s,\left( x_{n}\right) \right) \right\} \subset U\left(
s\right) .
\]
Let 
\[
W\stackrel{\text{def}}{=}\left\{ s\in \mathbf{R}^{\kappa }\mid \tilde{\alpha}
\left( s,\left( x\right) \right) \text{ and }h^{\left[ x_{n},x\right] }\circ 
\tilde{\alpha}\left( s,\left( x_{n}\right) \right) \text{ are in the same
component of }U\left( s\right) \right\} .
\]
By construction, both $W$ and $\mathbf{R}^{\kappa }-W$ are open in $\mathbf{R
}^{\kappa }$ since the components of $U\left( s\right) $ are homeomorphic to 
$U$. Since $h^{\left[ x_{n},x\right] }\left( x_{n}\right) =x$, $\mathbf{0}
\in W$ and since $\mathbf{R}^{\kappa }$ is connected, $W=\mathbf{R}^{\kappa }
$. Thus, with $C\left( s\right) $ denoting the component of $U\left(
s\right) $ containing both $\tilde{\alpha}\left( s,\left( x\right) \right) $
and $h^{\left[ x_{n},x\right] }\circ \tilde{\alpha}\left( s,\left(
x_{n}\right) \right) $, and with 
\[
f_{s}\stackrel{\text{def}}{=}C\left( s\right) \cap p^{-1}\left( \alpha
\left( s,p\left( x_{n}\right) \right) \right) 
\]
we have that 
\[
\left\{ \tilde{\alpha}\left( s,x\right) ,h^{\left[ x_{n},x\right] }\circ 
\tilde{\alpha}\left( s,\left( x_{n}\right) \right) \right\} \subset \tilde{
\alpha}\left( U\times \left\{ f_{s}\right\} \right) ,
\]
implying that $d\left( \tilde{\alpha}\left( s,x\right) ,h^{\left[ x_{n},x
\right] }\circ \tilde{\alpha}\left( s,\left( x_{n}\right) \right) \right)
<\varepsilon /3$ by our initial choice of $U$. The choice of $h^{\left[
x_{n},x\right] }$ yields that $d\left( \tilde{\alpha}\left( s,x\right) ,
\tilde{\alpha}\left( s,\left( x_{n}\right) \right) \right) <5\varepsilon /6$
. Since $s$ was any point of $\mathbf{R}^{\kappa }$ , we may finally
conclude that 
\[
\sup\limits_{s\in \mathbf{R}^{\kappa }}d\left( \tilde{\alpha}_{x_{n}}\left(
s\right) ,\tilde{\alpha}_{x}\left( s\right) \right) <\varepsilon \text{ for
all }n\geq N.
\]
And so $h$ is a homeomorphism. It follows easily that $\tilde{\alpha}$ is
uniformly Lyapunov stable in the strongest sense: for any given $\varepsilon
>0$ there is a $\delta >0$ so that 
\[
d\left( x,y\right) <\delta \Rightarrow d\left( \tilde{\alpha}\left(
s,x\right) ,\tilde{\alpha}\left( s,y\right) \right)  <\varepsilon  \text{ for all }s\in 
\mathbf{R}^{\kappa }.
\]

Fixing a point $e\in X$, the density of the orbit of $e$ means that any $
a,b\in X$ may be represented in the form $a=\lim_{i}\left\{ \tilde{\alpha}
\left( t_{i}^{a},e\right) \right\} $ and $b=\lim_{i}\left\{ \tilde{\alpha}
\left( t_{i}^{b},e\right) \right\} $. The operations 
\[
-a=\lim_{i}\left\{ \tilde{\alpha}\left( -t_{i}^{a},e\right) \right\} \text{
and }a+b=\lim_{i}\left\{ \tilde{\alpha}\left( t_{i}^{a}+t_{i}^{b},e\right)
\right\} 
\]
then give $X$ a well-defined, abelian topological group structure compatible
with the original topology of $X$. The proof is essentially identical with
that found in \cite{NS} V, 8.16 and is therefore omitted.\hfill $\square $

It follows that such an $X$ is homeomorphic to either $\mathbf{T}^{\kappa }$
or the inverse limit of an inverse sequence of $\mathbf{T}^{\kappa }$ with
epimorphic bonding maps (see \cite{C}). And it follows directly that any
such $X$ is bihomogeneous: the homeomorphism 
\[
x\mapsto \left( a+b\right) -x
\]
switches $a$ and $b$.

\section{Exponents}

We now define the exponent group and explore its topological significance.
In what follows, $Hom\left( \mathbf{R}^{\kappa },\mathbf{R}^{\kappa }\right) $
denotes the group of continuous homomorphisms $\mathbf{R}^{\kappa }\rightarrow 
\mathbf{R}^{\kappa }$ with point-wise addition and $\left[ X;Y\right] $ denotes
the group of homotopy classes of maps $X\rightarrow Y$, and $\left[ f\right] 
$ denotes the homotopy class of $f:X\rightarrow Y$.

\begin{Definition}
\emph{Given a map }$f:W\rightarrow X$\emph{, we define }$\left\{ w_i\right\}
_{i=1}^\infty \subset W$ \emph{to be an }$f$\textit{--sequence} \emph{if} $
\left\{ f\left( w_i\right) \right\} _{i=1}^\infty $ \emph{converges in }$X$.
\end{Definition}

\begin{Definition}
\emph{Given a map }$f$\emph{\ of }$\mathbf{R}^\kappa $\emph{\ into a metric
space }$X$\emph{, the }\textit{exponent} \textit{group of}\emph{\ }$f$ 
\textit{,}\emph{\ denoted }$\mathcal{E}_f,$\emph{\ is } 
\[
\left\{ A\in Hom\left( \mathbf{R}^\kappa ,\mathbf{R}^\kappa \right) \mid \;\left\{
\pi ^\kappa \left( A\left( t_i\right) \right) \right\} _{i=1}^\infty \text{
\emph{converges in} }\mathbf{T}^\kappa \text{\emph{for all} }f\text{\textit{
\ --\emph{sequences} }}\left\{ t_i\right\} _{i=1}^\infty \right\} .
\]
\end{Definition}

\begin{Lemma}
$\mathcal{E}_f$ is a subgroup of $Hom\left( \mathbf{R}^\kappa ,\mathbf{R}^\kappa
\right) $.
\end{Lemma}

\noindent \textbf{Proof}: The zero homomorphism is trivially in $\mathcal{E}
_{f}$. If $A$ and $B$ are in $\mathcal{E}_{f}$ and if $\left\{ t_{i}\right\}
_{i=1}^{\infty }$ is an $f$\textit{--}sequence, then $\left\{ \pi ^{\kappa
}\left( A\left( t_{i}\right) \right) \right\} _{i=1}^{\infty }$ and $\left\{
\pi ^{\kappa }\left( B\left( t_{i}\right) \right) \right\} _{i=1}^{\infty }$
converge in $\mathbf{T}^{\kappa }$. And since $-$ is continuous on $
Hom\left( \mathbf{R}^{\kappa },\mathbf{R}^{\kappa }\right) \times Hom\left( \mathbf{R}
^{\kappa },\mathbf{R}^{\kappa }\right) $, this implies that $\left\{ \pi
^{\kappa }\left( \left( A-B\right) \left( t_{i}\right) \right) \right\}
_{i=1}^{\infty }$ converges in $\mathbf{T}^{\kappa }$.\hfill $\square $

\begin{Theorem}
Given a map $f$ of $\mathbf{R}^\kappa $ into a metric space $X$, with 
\[
\iota :\mathcal{E}_f\rightarrow \left[ \overline{f\left( \mathbf{R}^\kappa
\right) };\mathbf{T}^\kappa \right]
\]
given by 
\[
A\stackrel{\iota }{\mapsto }\left[ f_A\right] \text{, where }f_A:\overline{
f\left( \mathbf{R}^\kappa \right) }\rightarrow \mathbf{T}^\kappa \text{ sends }
\lim_i\left\{ f\left( t_i\right) \right\} \text{ to }\lim_i\left\{ \pi
^\kappa \left( A\left( t_i\right) \right) \right\}
\]
we have:
\end{Theorem}

\begin{enumerate}
\item  $\iota $ is a homomorphism (we give $\left[ \overline{f\left( \mathbf{R}
^{\kappa }\right) };\mathbf{T}^{\kappa }\right] $ the group operation
induced by point-wise addition of maps).

\item  If $\iota \left( A\right) =\iota \left( B\right) $ and if $\mathcal{C}
$ is any contractible topological subspace of $\mathbf{R}^\kappa $ satisfying
the condition that the map $\left( A-B\right) _{\mathcal{C}}:\mathcal{C}
\rightarrow \left( A-B\right) \left( \mathcal{C}\right) $; $t\mapsto \left(
A-B\right) \left( t\right) $ has a continuous inverse, then for any $f$
--sequence $\left\{ t_i\right\} _{i=1}^\infty \subset \mathcal{C} $, the
sequence $\left\{ \left( A-B\right) \left( t_i\right) \right\} _{i=1}^\infty 
$ converges in $\mathbf{R}^\kappa $.

\item  If $\overline{f\left( \mathbf{R}^\kappa \right) }$ is compact, then $
\iota $ is an embedding and $\mathcal{E}_f$ is countable.

\item  If $\iota \left( A\right) =\iota \left( B\right) $ and if $\kappa
=n<\infty $ and if there is an unbounded $f$--sequence $\left\{ t_i\right\}
_{i=1}^\infty $, then $A-B$ is not invertible.
\end{enumerate}

\noindent \textbf{Proof}: $(1)$ Let $A\in \mathcal{E}_f$ be given. If $
\lim\limits_i\left\{ f\left( s_i\right) \right\} =\lim\limits_i\left\{
f\left( t_i\right) \right\} =x$ are two representations of a point in $
\overline{f\left( \mathbf{R}^\kappa \right) }$, then $\lim\limits_i\left\{ \pi
^\kappa \left( A\left( s_i\right) \right) \right\} =\sigma $ and $
\lim\limits_i\left\{ \pi ^\kappa \left( A\left( t_i\right) \right) \right\}
=\tau $ both exist by the definition of $\mathcal{E}_f$. Then for $i\in
\left\{ 1,2,...\right\} $ if we define $u_{2i-1}=s_i$ and $u_{2i}=t_i$, we
have that $\lim\limits_i\left\{ f\left( u_i\right) \right\} =x$, implying
that $\lim\limits_i\left\{ \pi ^\kappa \left( A\left( u_i\right) \right)
\right\} $ exists, which is only possible if $\sigma =\tau $. Thus, $f_A$ is
a well--defined function. To see that $f_A$ is continuous, consider a
convergent sequence $\left\{ \lim\limits_i\left\{ f\left( t_i^j\right)
\right\} \right\} _{j=1}^\infty =\left\{ x^j\right\} _{j=1}^\infty
\rightarrow x$. Then with $d$ denoting a metric for $X$ and with $d_{\kappa}$ 
denoting a metric for $\mathbf{T}^{\kappa}$, for each $j\in \left\{
1,2,...\right\} $ choose $i_j$ so that: 
\[
d\left( f\left( t_{i_j}^j\right) ,x^j\right) <\frac 1j\text{ and }d_{\kappa}
\left( f_A\left( x^j\right) ,\pi ^\kappa \left( A\left(
t_{i_j}^j\right) \right) \right) <\frac 1j.
\]
Then $\lim\limits_j\left\{ f\left( t_{i_j}^j\right) \right\} _{j=1}^\infty
=x $ and 
\[
f_A\left( x\right) =\lim\limits_j\left\{ \pi ^\kappa \left( A\left(
t_{i_j}^j\right) \right) \right\} _{j=1}^\infty =\lim\limits_j\left\{
f_A\left( x^j\right) \right\} _{j=1}^\infty ,
\]
demonstrating that $f_A$ is continuous. And if $A,B\in \mathcal{E}_f$ we
have for $x=\lim\limits_i\left\{ f\left( t_i\right) \right\} $ 
\[
f_{A+B}\left( x\right) =\lim\limits_i\left\{ \pi ^\kappa \left( \left(
A+B\right) \left( t_i\right) \right) \right\} =\lim\limits_i\left\{ \pi
^\kappa \left( A\left( t_i\right) \right) \right\} +\lim\limits_i\left\{ \pi
^\kappa \left( B\left( t_i\right) \right) \right\} =f_A\left( x\right)
+f_B\left( x\right) \text{,}
\]
demonstrating that $\iota $ is a homomorphism.

$(2)$ Since $\left[ f_A\right] =\left[ f_B\right] $, we have that $\left[
f_{\left( A-B\right) }\right] =\left[ \text{constant map}\right] $, and so
there is a map $g:\left( \overline{f\left( \mathbf{R}^\kappa \right) },f\left( 
\mathbf{0}\right) \right) \rightarrow \left( \mathbf{R}^\kappa ,\mathbf{0}
\right) $ lifting $f_{\left( A-B\right) }$ making the following diagram
commute 
\[
\begin{array}{ccccc}
&  & \mathbf{R}^\kappa &  &  \\ 
& ^g\nearrow &  & \searrow ^{\pi ^\kappa } &  \\ 
\overline{f\left( \mathbf{R}^\kappa \right) } &  & \stackrel{f_{\left(
A-B\right) }}{\longrightarrow } &  & \mathbf{T}^\kappa
\end{array}
.
\]
Then for any $t\in \left( A-B\right) \left( \mathcal{C}\right) $ we have 
\[
f_{\left( A-B\right) }\circ f\circ \left( A-B\right) _{\mathcal{C}
}^{-1}\left( t\right) =\pi ^\kappa \left( \left( A-B\right) \left(
A-B\right) _{\mathcal{C}}^{-1}\left( t\right) \right) =\pi ^\kappa \left(
t\right) ,
\]
and so we are led to the following commutative diagram: 
\[
\begin{array}{ccccccc}
\left( A-B\right) \left( \mathcal{C}\right) & \stackrel{\left( A-B\right) _{
\mathcal{C}}^{-1}}{\rightarrow } & \mathcal{C} & \stackrel{f}{\rightarrow }
& f\left( \mathcal{C}\right) & \stackrel{g}{\rightarrow } & \mathbf{R}^\kappa
\\ 
\pi ^\kappa \downarrow &  &  &  &  & \searrow f_{\left( A-B\right) } & 
\downarrow \pi ^\kappa \\ 
\mathbf{T}^\kappa &  & \stackrel{id}{\longrightarrow } &  &  &  & \mathbf{T}
^\kappa
\end{array}
\text{.}
\]
And since, $g\circ f\circ \left( A-B\right) _{\mathcal{C}}^{-1}$ and $id_{
\mathbf{R}^\kappa }$ both map $\mathbf{0}$ to $\,\mathbf{0}$, we have that $
g\circ f\circ \left( A-B\right) _{\mathcal{C}}^{-1}=id_{\mathbf{R}^\kappa }$
since both provide a lift of $id_{\mathbf{T}^1}\circ \pi ^\kappa \mid
_{\left( A-B\right) \left( \mathcal{C}\right) }$ and such a lift is uniquely
determined. And so if $\left\{ t_i\right\} _{i=1}^\infty \subset \mathcal{C}$
is an $f$--sequence with $\lim\limits_i\left\{ f\left( t_i\right) \right\}
=x\in \overline{f\left( \mathbf{R}^\kappa \right) }$, we must have 
\begin{eqnarray*}
g\left( x\right) &=&g\left( \lim_i\left\{ f\circ \left( A-B\right) _{
\mathcal{C}}^{-1}\left( \left( A-B\right) _{\mathcal{C}}\left( t_i\right)
\right) \right\} \right) \\
&=&\lim_i\left\{ g\circ f\circ \left( A-B\right) _{\mathcal{C}}^{-1}\left(
\left( A-B\right) _{\mathcal{C}}\left( t_i\right) \right) \right\}
=\lim_i\left\{ \left( A-B\right) _{\mathcal{C}}\left( t_i\right) \right\} 
\text{ in }\mathbf{R}^\kappa
\end{eqnarray*}
by the continuity of $g$.

$(3)$ Suppose then that $\overline{f\left( \mathbf{R}^\kappa \right) }$ is
compact and that $\iota \left( A\right) =\left[ f_A\right] =\left[ \text{
constant map}\right] $ and that $v\in \mathbf{R}^\kappa -\ker A$. Then with $
\mathcal{C}$ denoting the vector subspace $\mathbf{R\cdot }v\subset \mathbf{R}
^\kappa $, we have that $A_{\mathcal{C}}=\left( A-0\right) _{\mathcal{C}}$
is an isomorphism onto $\mathbf{R\cdot }A\left( v\right) $. And since $
\overline{f\left( \mathbf{R}^\kappa \right) }$ is compact, there is a
subsequence $\left\{ f\left( n_iv\right) \right\} _{i=1}^\infty $ of the
sequence $\left\{ f\left( nv\right) \right\} _{n=1}^\infty $ which converges
to some $x\in \overline{f\left( \mathbf{R}^\kappa \right) }$. But then by the
above, we must have that $\left\{ n_iA\left( v\right) \right\} _{i=1}^\infty 
$ converges in $\mathbf{R}^\kappa $. And so if $w_\ell $ is any non-zero
component of $A\left( v\right) =\left\langle w_i\right\rangle _{i=1}^\kappa $
, we would then have that $\left\{ n_iw_\ell \right\} _{i=1}^\infty $
converges in $\mathbf{R}$, which is impossible since $\left\{ n_i\right\}
_{i=1}^\infty $ is unbounded. Thus, we must have $\ker A=\mathbf{R}^\kappa $
and $A$ is the zero map. That $\mathcal{E}_f$ is countable then follows from
the fact that $\left[ \overline{f\left( \mathbf{R}^\kappa \right) };\mathbf{T}
^\kappa \right] $ is countable.

$\left( 4\right) $ If $\kappa <\infty $ and if $\left\{ t_{i}\right\}
_{i=1}^{\infty }$ is an unbounded $f$--sequence, then $\left\{ \left(
A-B\right) \left( t_{i}\right) \right\} _{i=1}^{\infty }$ would be unbounded
if $\left( A-B\right) $ were invertible.\hfill $\square $

We know that when $f$ has an image which is not compact, $\mathcal{E}_f$ may
not be countable and $\iota $ may not be an embedding. For example, for $f: 
\mathbf{R}^2\rightarrow \mathbf{R}\times \mathbf{T}^1;$ $\left\langle
t_{1,}t_2\right\rangle \mapsto \left\langle t_1,\pi ^1\left( t_2\right)
\right\rangle $, making the identification of $Hom\left( \mathbf{R}^2,\mathbf{R}
^2\right) $ with $Mat\left( 2,\mathbf{R}\right) $ (the group of $2\times 2$
matrices with real entries under component--wise addition), we have 
\[
\left\{ \left( 
\begin{array}{ll}
c & 0 \\ 
0 & z
\end{array}
\right) \in Mat\left( 2,\mathbf{R}\right) \mid c\in \mathbf{R}\text{ and }z\in 
\mathbf{Z}\right\} \subset \mathcal{E}_f
\]
and $\iota \left( \left( 
\begin{array}{ll}
c & 0 \\ 
0 & z
\end{array}
\right) \right) =\iota \left( \left( 
\begin{array}{ll}
c^{\prime } & 0 \\ 
0 & z
\end{array}
\right) \right) $ for any $c$ and $c^{\prime }$ in $\mathbf{R}$. And since
there are unbounded $f$--sequences, we know that $\left( 
\begin{array}{ll}
c & 0 \\ 
0 & z
\end{array}
\right) -\left( 
\begin{array}{ll}
c^{\prime } & 0 \\ 
0 & z
\end{array}
\right) $ is not invertible.

\begin{Theorem}
\label{Map}If $f:\mathbf{R}^\kappa \mathbf{\rightarrow }X$ has exponent group $
\mathcal{E}_f=\left\{ A_\lambda \mid \;\lambda \in \Lambda \right\} $, we
have the map 
\[
h_f:\overline{f\left( \mathbf{R}^\kappa \right) }\rightarrow
\prod\limits_{\lambda \in \Lambda }\mathbf{T}^\kappa ;\text{ }x\mapsto
\left\langle f_{A_\lambda }\left( x\right) \right\rangle _{\lambda \in
\Lambda },
\]
and $h_f\circ f$ is a homomorphism of $\mathbf{R}^\kappa $ into $
\prod\limits_{\lambda \in \Lambda }\mathbf{T}^\kappa $. And if $\overline{
f\left( \mathbf{R}^\kappa \right) }$ is compact, then $h_f\left( \overline{
f\left( \mathbf{R}^\kappa \right) }\right) $ is a subgroup of $
\prod\limits_{\lambda \in \Lambda }\mathbf{T}^\kappa $.
\end{Theorem}

\noindent \textbf{Proof}: Since each $f_{A_{\lambda }}$ is continuous, we
know that $h_{f}$ is also continuous. When we give $\prod\limits_{\lambda
\in \Lambda }\mathbf{T}^{\kappa }$ the group structure of component--wise
addition, it is a topological group with the product (Tychonoff) topology
and 
\begin{eqnarray*}
h_{f}\circ f\left( s+t\right)  &=&\left\langle \pi ^{\kappa }\circ
A_{\lambda }\left( s+t\right) \right\rangle _{\lambda \in \Lambda
}=\left\langle \pi ^{\kappa }\circ A_{\lambda }\left( s\right) \right\rangle
_{\lambda \in \Lambda }+\left\langle \pi ^{\kappa }\circ A_{\lambda }\left(
t\right) \right\rangle _{\lambda \in \Lambda } \\
&=&h_{f}\circ f\left( s\right) +h_{f}\circ f\left( t\right) ,
\end{eqnarray*}
and so $h_{f}\circ f$ is a homomorphism. When $\overline{f\left( \mathbf{R}
^{\kappa }\right) }$ is compact, we have that $h_{f}\left( \overline{f\left( 
\mathbf{R}^{\kappa }\right) }\right) =\overline{h_{f}\circ f\left( \mathbf{R}
^{\kappa }\right) }$ (recall that $\Lambda $ is countable when $\overline{
f\left( \mathbf{\ R}^{\kappa }\right) }$ is compact). We now know that $
h_{f}\circ f\left( \mathbf{R}^{\kappa }\right) $ is a subgroup of $
\prod\limits_{\lambda \in \Lambda }\mathbf{T}^{\kappa }$, and so $\overline{
h_{f}\circ f\left( \mathbf{R}^{\kappa }\right) }$ is also a subgroup.\hfill $
\square $

We are now in a position to give a generalization of almost periodic maps $
\mathbf{R\rightarrow }X$ to almost periodic maps $\mathbf{R}^{\kappa }\
\rightarrow X$.

\begin{Definition}
\emph{A map} $f:\mathbf{R}^\kappa \mathbf{\rightarrow }X$ \emph{is} \textit{almost
periodic} \emph{if} $\mathcal{E}_f$ \emph{is countable and if the following
condition} \emph{is satisfied}: 
\[
\left( \ast \right) \;\left[ \left\{ t_i\right\} _{i=1}^\infty \text{\emph{\
is an} }f\text{--\emph{sequence}}\right] \Leftrightarrow \left[ \left\{ \pi
^\kappa \left( A\left( t_i\right) \right) \right\} _{i=1}^\infty \text{\emph{
\ converges for all} }A\in \mathcal{E}_f\right] .
\]
\end{Definition}

It is known that any classically defined almost periodic function satisfies $
\left( \ast \right) $ when $\kappa =1$ and $X$ is complete, see \cite{Cl}
. And below we shall see that any map satisfying $\left( \ast \right) $ is
indeed almost periodic in the classical sense when $\kappa =1$ and similar
techniques may be used to show that generally the above definition coincides
with Bochner's \cite{B}.

\begin{Theorem}
If $f$ is almost periodic, then $\overline{f\left( \mathbf{R}^\kappa \right) }$
is compact, and $h_f$ as in Theorem \ref{Map} is a homeomorphism of $
\overline{f\left( \mathbf{R}^\kappa \right) }$ onto the subgroup $h_f\left( 
\overline{f\left( \mathbf{R}^\kappa \right) }\right) $ of $\prod\limits_{n=1}^
\infty \mathbf{T}^\kappa $.
\end{Theorem}

\noindent \textbf{Proof}: Let $\left\{ x^{j}\right\} _{j=1}^{\infty }$ be a
sequence in $\overline{f\left( \mathbf{R}^{\kappa }\right) }$. Choose a
subsequence of $\left\{ h_{f}\left( x^{j}\right) \right\} _{j=1}^{\infty }$
which converges in $\prod\limits_{n=1}^{\infty }\mathbf{T}^{\kappa }$ and
which we label $\left\{ h_{f}\left( y^{j}\right) \right\} _{j=1}^{\infty }$
for convenience. With $y^{j}=\lim\limits_{i}\left\{ f\left( t_{i}^{j}\right)
\right\} _{i=1}^{\infty }$ and with $\mathcal{E}_{f}=\left\{
A_{1},A_{2},...\right\} ,$ for each $j\in \left\{ 1,2,...\right\} $ choose $
t_{i_{j}}^{j}$ such that 
\[
d\left( f\left( t_{i_{j}}^{j}\right) ,y^{j}\right) <\frac{1}{j}\text{ and }
\max \left\{ d_{\mathbf{T}^{\kappa }}\left( \pi ^{\kappa }\left( A_{n}\left(
t_{i_{j}}^{j}\right) \right) ,f_{A_{n}}\left( y^{j}\right) \right) \mid n\in
\left\{ 1,...,j\right\} \right\} <\frac{1}{j}\text{.}
\]
Since we have that $\left\{ f_{A_{n}}\left( y^{j}\right) \right\}
_{j=1}^{\infty }$ converges for each $n\in \left\{ 1,2,...\right\} $, we
must then have that $\left\{ \pi ^{\kappa }\left( A\left(
t_{i_{j}}^{j}\right) \right) \right\} _{j=1}^{\infty }$ converges for all $
A\in \mathcal{E}_{f}$, and so $\left\{ t_{i_{j}}^{j}\right\} _{j=1}^{\infty }
$ is an $f$--sequence by $\left( \ast \right) $, say $\left\{ f\left(
t_{i_{j}}^{j}\right) \right\} _{j=1}^{\infty }\rightarrow y\in X$. And so we
also have $\left\{ f\left( y^{j}\right) \right\} _{j=1}^{\infty }\rightarrow
y$, demonstrating that $\overline{f\left( \mathbf{R}^{\kappa }\right) }$ is
compact. Now suppose that $h_{f}\left( x\right) =h_{f}\left( y\right) $ with 
$x=\lim\limits_{i}\left\{ f\left( s_{i}\right) \right\} $ and $
y=\lim\limits_{i}\left\{ f\left( t_{i}\right) \right\} $. We then have that
for each $n\in \left\{ 1,2,...\right\} $ 
\[
f_{A_{n}}\left( x\right) =\lim_{i}\left\{ \pi ^{\kappa }\left( A_{n}\left(
s_{i}\right) \right) \right\} _{i=1}^{\infty }=\lim_{i}\left\{ \pi ^{\kappa
}\left( A_{n}\left( t_{i}\right) \right) \right\} _{i=1}^{\infty
}=f_{A_{n}}\left( y\right) ,
\]
and so with $u_{2i-1}=s_{i}$ and $u_{2i}=t_{i}$ for $i\in \left\{
1,2,...\right\} ,$ we have that $\lim\limits_{i}\left\{ \pi ^{\kappa }\left(
A_{n}\left( u_{i}\right) \right) \right\} _{i=1}^{\infty }=f_{A_{n}}\left(
x\right) $ for each $n\in \left\{ 1,2,...\right\} .$ Then by $\left( \ast
\right) $ we have that $\left\{ f\left( u_{i}\right) \right\} _{i=1}^{\infty
}$ converges in $X$, which is only possible if $x=y$. Thus, $h_{f}$ is
one--to--one and so is a homeomorphism onto its image.\hfill $\square $

When $f$ is almost periodic, there is a topological isomorphism $\mathfrak{i}_{f}
$ from $h_{f}\left( \overline{f\left( \mathbf{R}^{\kappa }\right) }\right) $
onto some $\lambda $--solenoid $\sum_{\overline{M}}$ $\left( \lambda \leq
\infty \right) $ since $h_{f}\left( \overline{f\left( \mathbf{R}^{\kappa
}\right) }\right) $ is a compact connected abelian group; see \cite{P}
Thm 68, \cite{NS} V,8.16 \cite{C} and \cite{Cl}. Then, using the notation of 
\cite{C}, we have a map $\mathfrak{h}:\left( \mathbf{R}^{\kappa },\mathbf{0}\right)
\rightarrow \left( \mathbf{R}^{\lambda },\mathbf{0}\right) $ making the
following diagram commute: 
\[
\begin{array}{ccc}
&  & \mathbf{R}^{\lambda } \\ 
& \mathfrak{h}\nearrow  & \downarrow \pi _{\overline{M}} \\ 
\mathbf{R}^{\kappa } & \stackrel{\mathfrak{i}_{f}\circ h_{f}\circ f}{
\longrightarrow } & \sum_{\overline{M}}
\end{array}
\]
since $\pi _{\overline{M}}$ is a fibration with unique path lifting. Then
for $s,t\in \mathbf{R}^{\kappa }$, since $\mathfrak{i}_{f}\circ h_{f}\circ f$ and $
\pi _{\overline{M}}$ are homomorphisms, we have 
\begin{eqnarray*}
\pi _{\overline{M}}\left( \mathfrak{h}\left( s\right) +\mathfrak{h}\left( t\right)
\right)  &=&\pi _{\overline{M}}\circ \mathfrak{h}\left( s\right) +\pi _{
\overline{M}}\circ \mathfrak{h}\left( t\right) =\mathfrak{i}_{f}\circ h_{f}\circ
f\left( s\right) +\mathfrak{i}_{f}\circ h_{f}\circ f\left( t\right)  \\
&=&\mathfrak{i}_{f}\circ h_{f}\circ f\left( s+t\right) =\pi _{\overline{M}}\circ 
\mathfrak{h}\left( s+t\right) ,
\end{eqnarray*}
and so $\mathfrak{h}\left( s\right) +\mathfrak{h}\left( t\right) -\mathfrak{h}\left(
s+t\right) \in \ker \pi _{\overline{M}}$, implying that $\mathfrak{h}\left(
s\right) +\mathfrak{h}\left( t\right) =\mathfrak{h}\left( s+t\right) $ for all $
s,t\in \mathbf{R}^{\kappa }$, see \cite{C} 3.4. Hence, $\mathfrak{h}\left( \mathbf{R}
^{\kappa }\right) $ is a subgroup of $\mathbf{R}^{\lambda }$. And so we may
think of $\pi _{\overline{M}}\circ \mathfrak{h}\left( \mathbf{R}^{\kappa }\right) $
as a ``linear subspace'' of $\sum_{\overline{M}}$.

\begin{Theorem}
If $f:\mathbf{R}^\kappa \rightarrow X$ is almost periodic, then $f$ may be
extended to a continuous group action of $\left( \mathbf{R}^\kappa ,+\right) $
on all of $\overline{f\left( \mathbf{R}^\kappa \right) }:$ 
\[
\alpha _f:\mathbf{R}^\kappa \times \overline{f\left( \mathbf{R}^\kappa \right) }
\rightarrow \overline{f\left( \mathbf{R}^\kappa \right) },\text{ where }\alpha
_f\left( t\mathbf{,}f\left( 0\right) \right) =f\left( t\right) .
\]
\end{Theorem}

\noindent \textbf{Proof}: We have the group action 
\[
\phi _{f}:\mathbf{R}^{\kappa }\times \sum\nolimits_{\overline{M}}\rightarrow
\sum\nolimits_{\overline{M}};\;\left( t,x\right) \mapsto \pi _{\overline{M}
}\circ \mathfrak{h}\left( t\right) +x.
\]
For $x\in \sum\nolimits_{\overline{M}}$, we have $\phi _{f}\left( \mathbf{0,}
x\right) =\pi _{\overline{M}}\circ \mathfrak{h}\left( \mathbf{0}\right) +x=e_{
\overline{M}}+x=x$ and for $s,t\in \mathbf{R}^{\kappa }$ we have 
\begin{eqnarray*}
\phi _{f}\left( s+t\mathbf{,}x\right)  &=&\pi _{\overline{M}}\circ \mathfrak{h}
\left( s+t\right) +x=\pi _{\overline{M}}\circ \mathfrak{h}\left( s\right) +\pi _{
\overline{M}}\circ \mathfrak{h}\left( t\right) +x \\
&=&\phi _{f}\left( s\mathbf{,}\phi _{f}\left( t,x\right) \right) ,
\end{eqnarray*}
and so $\phi _{f}$ is indeed a group action. We then define 
\[
\alpha _{f}\left( s,x\right) \stackrel{def}{=}\left( \mathfrak{i}_{f}\circ
h_{f}\right) ^{-1}\left( \phi _{f}\left( s,\mathfrak{i}_{f}\circ h_{f}\left(
x\right) \right) \right) ,
\]
i.e., the action conjugate via $\mathfrak{i}_{f}\circ h_{f}$ to $\phi _{f}$. It
then follows directly that $\alpha _{f}$ is a group action, and 
\begin{eqnarray*}
\alpha _{f}\left( t\mathbf{,}f\left( 0\right) \right)  &=&\left( \mathfrak{i}
_{f}\circ h_{f}\right) ^{-1}\left( \phi _{f}\left( t,\mathfrak{i}_{f}\circ
h_{f}\left( f\left( 0\right) \right) \right) \right) =\left( \mathfrak{i}
_{f}\circ h_{f}\right) ^{-1}\left( \pi _{\overline{M}}\circ \mathfrak{h}\left(
t\right) +e_{\overline{M}}\right)  \\
&=&\left( \mathfrak{i}_{f}\circ h_{f}\right) ^{-1}\left( \mathfrak{i}_{f}\circ
h_{f}\circ f\left( t\right) \right) =f\left( t\right). \square 
\end{eqnarray*}

And when $\kappa =1$, $\alpha _{f}$ is a flow on $\overline{f\left( \mathbf{R}
^{\kappa }\right) }$ and $\mathfrak{i}_{f}\circ h_{f}$ provides an equivalence
of $\alpha _{f}$ and an almost periodic linear flow as described in \cite{C}
, implying that $f$ is itself almost periodic in the classical sense since
it is an orbit of this flow.

Generally, it is now clear that all $\phi _{f}$--orbits are translates of $
\pi _{\overline{M}}\circ \mathfrak{h}\left( \mathbf{R}^{\kappa }\right) $, and so
we consider the decomposition of $\sum\nolimits_{\overline{M}}$ determined
by $\phi _{f}$ as a decomposition into ``linear subspaces.'' Thus, $f$
determines a group action on $\overline{f\left( \mathbf{R}^{\kappa }\right) }$
which is equivalent to a ``linear foliation'' of $\sum\nolimits_{\overline{M}
}$. It then follows that if $f:\mathbf{R}^{n}\rightarrow M$ is a smooth almost
periodic map onto a leaf of an $n$--dimensional foliation $\mathcal{F}$ of
the manifold $M$ that $\mathfrak{i}_{f}\circ h_{f}$ then provides an equivalence
of $\mathcal{F}\mid _{\overline{f\left( \mathbf{R}^{n}\right) }}$ with the
linear foliation $\phi _{f}$ since the leaves of $\mathcal{F}\mid _{
\overline{f\left( \mathbf{R}^{n}\right) }}$ will coincide with the orbits of $
\alpha _{f}$ due to the density of $f\left( \mathbf{R}^{n}\right) $ in $
\overline{f\left( \mathbf{R}^{n}\right) }$.

Given $\mathbf{v\in R}^{\kappa }$, any map $f:\mathbf{R}^{\kappa
}\rightarrow X$ induces the map 
\[
\mathbf{R\approx \,}t\cdot \mathbf{v\hookrightarrow R}^{\kappa }\stackrel{f}{
\rightarrow }X
\]
and the one--dimensional exponent group of this map $\mathbf{R\rightarrow }X$
is frequently sufficient to determine as much about $\check{H}^{1}\left( 
\overline{f\left( \mathbf{R}^{\kappa }\right) }\right) \approx \left[ 
\overline{f\left( \mathbf{R}^{\kappa }\right) };\mathbf{T}^{1}\right] $ as $
\mathcal{E}_{f}$ reveals. But when $f$ extends to a group action,  $f_{A}$
will be a fiber bundle projection, and some of the bundle structure is lost
in considering only the one--dimensional maps. Thus, the higher dimensional
groups have content not captured in the one-dimensional exponent group.

\end{document}